\documentclass[a4paper,12pt,reqno]{amsart}
\usepackage[T1]{fontenc}
\usepackage[utf8]{inputenc}
\usepackage[libertinus,vvarbb]{newtx}

\usepackage[margin=1in]{geometry}
\usepackage{enumitem}
\usepackage{graphicx}
\usepackage[dvipsnames]{xcolor}
\usepackage{tikz}
\usepackage{xparse}
\usepackage{xurl}
\usepackage[hidelinks]{hyperref}
\usepackage{tikz}
\usepackage{xcolor}
\usepackage{cite}
\usepackage{changes} 

\definecolor{MyBlue}{RGB}{200, 225, 255}
\definecolor{MyGreen}{RGB}{200, 245, 200}
\definecolor{MyOrange}{RGB}{255, 225, 180}

\usetikzlibrary{decorations.markings,decorations.pathmorphing,nfold}
\makeatletter
\tikzset{
 side by side path/.style args={#1:#2}{
  postaction={path only,draw=#1,offset=+.5\pgflinewidth},
  postaction={path only,draw=#2,offset=+-.5\pgflinewidth}
 },
 side by side/.style={path only,side by side path={#1}},
 offset/.code=\tikz@addoption{%
  \pgfgetpath\tikz@temp%
  \pgfsetpath\pgfutil@empty%
  \pgfoffsetpath\tikz@temp{#1}%
 }
}
\makeatother

\numberwithin{equation}{section}

\newtheorem{theorem}{Theorem}[section]
\newtheorem{lemma}[theorem]{Lemma}
\newtheorem{corollary}[theorem]{Corollary}
\newtheorem{proposition}[theorem]{Proposition}

\theoremstyle{definition}

\newtheorem{example}[theorem]{Example}
\newtheorem{remark}[theorem]{Remark}

\DeclareMathOperator{\re}{Re}
\DeclareMathOperator{\im}{Im}

\newcommand{\ind}{\mathbb{1}}

\newcommand{\pr}{\mathbb{P}}

\newcommand{\C}{\mathbb{C}}
\newcommand{\R}{\mathbb{R}}
\newcommand{\E}{\mathbb{E}}
\newcommand{\Z}{\mathbb{Z}}
\newcommand{\N}{\mathbb{N}}
\newcommand{\G}{\mathscr{G}}

\newcommand{\ph}{\varphi}

\mathchardef\mathhyphen="2D

\newcommand{\cm}{\mathscr{C}\mspace{-4mu}\mathscr{M}}
\newcommand{\am}{\mathscr{A}\mspace{-4mu}\mathscr{M}}
\newcommand{\amcm}{\am\mspace{-2mu}\mathhyphen\mspace{-1mu}\cm}

\renewcommand{\leq}{\leqslant}
\renewcommand{\geq}{\geqslant}

\newcommand{\lv}{\lvert}
\newcommand{\rv}{\rvert}

\newcommand{\floor}[1]{\left\lfloor #1 \right\rfloor}

\NewDocumentCommand{\formula}{ssom}{%
 \IfBooleanTF{#1}{%
  \IfBooleanTF{#2}{%
   \IfValueTF{#3}%
    {\begin{align}\label{#3}\begin{gathered}#4\end{gathered}\end{align}}%
    {\begin{gather}#4\end{gather}}%
  }{%
   \IfValueTF{#3}%
    {\begin{align}\label{#3}\begin{aligned}#4\end{aligned}\end{align}}%
    {\begin{gather*}#4\end{gather*}}%
  }%
 }{%
  \IfValueTF{#3}%
   {\begin{align}\label{#3}#4\end{align}}%
   {\begin{align*}#4\end{align*}}%
 }%
}
\begin{document}

\title[Characterisation of first passage locations]{Characterisation of first passage locations for two-dimensional lattice random walks}
\author{Jacek Wszoła}
\thanks{Work supported by the Polish National Science Centre (NCN) grant no.\@ 2023/49/B/ST1/04303}
\address{\normalfont Jacek Wszoła \\ Department of Analysis and Stochastic Processes \\ Wrocław University of Science and Technology \\ Wybrzeże Wyspiańskiego 27 \\ 50-370 Wrocław, Poland}
\email{jacek.wszola@pwr.edu.pl}
\date{\today}
\begin{abstract}
Consider a diagonal random walk $(X_n,Y_n)$ on $\Z^2$ with transition probabilities depending only on the second coordinate. We characterise the distributions of first passage locations, $X(\tau_a)$, where $\tau_a$ is the first time $Y_n$ reaches level $a$. We obtain a necessary and sufficient condition for a given distribution to be realised as a rescaled first passage location from some starting level. The proof extends Bondesson's method for one-dimensional first passage times. It combines a transformation of generating functions based on Pick functions with a construction using barriers. We also study maximal representations of generating functions of first passage locations and show that analogous converse results for standard and honeycomb random walks require additional restrictions.
\end{abstract}

\maketitle

\section{Introduction}

\subsection{Main results}

Let $(X_n)$ be a discrete-time random walk on $\Z$ with jumps $\pm1$. Assume that the probabilities of jumps depend only on the current position of the walk. By $\tau_0 = \min\{n \geq 0: X_n=0\}$ we denote the first passage time to $0$. The classical result of Bondesson \cite{bondesson} states that if the random walk starts from $x \in \N$, then $\tau_0$ has the same distribution as $x+2(X+Y)$, where $X$ has a distribution given by a mixture of geometric distributions and $Y$, independent of $X$, is a sum of at most $\floor{(x-1)/2}$ geometrically distributed random variables. Bondesson also proved the converse: for every distribution described above, there is a random walk such that $\tau_0$ restricted to $\{\tau_0 < \infty\}$ has the same distribution as $x+2(X+Y)$.

In this paper we study the two-dimensional variant of this setting. We consider a two-dimensional random walk $(X_n, Y_n)$ on $\Z^2$. Our process moves diagonally, with jumps $(\pm 1, \pm1)$. We additionally assume that the transition probabilities depend only on the current level $Y_n$ of the process. If $\tau_0 = \min\{n \geq 0: Y_n = 0\}$ is the first passage time to a level $0$, the random variable $X(\tau_0)$ is called the first passage location. Note that if the walk moves only to the right and $X_0=0$, then $X(\tau_0) = \tau_0$ and Bondesson's setting is recovered. 

In \cite{fpl}, a convolution representation of the first passage location distribution was obtained. More precisely, it was proved that if the random walk starts from $(0,y)$ and $y\geq 0$, then $(X(\tau_a)-y)/2$ restricted to $\{\tau_0<\infty\}$ has a distribution given by a convolution of geometric sequences, reflected geometric sequences, reflected two-point sequences and an $\amcm$ (absolutely monotone-then-completely monotone) sequence (see Theorem~1.1 in \cite{fpl}). The class of $\amcm$ sequences can be considered as a two-sided analogue of mixtures of geometric distributions (for details, see Section \ref{sec:amcm}). This raises the converse question of which distributions with this convolution structure can be realised by a diagonal random walk.

In this paper, we extend Bondesson's method to the two-dimensional setting and prove the following characterisation theorem.

\begin{theorem}
\label{thm:main}
    Let
    \formula[eq:main:F]{
    F(x) = G(x) \prod_{k=1}^{A} \frac{1-q_k}{1-q_k x} \prod_{k=1}^{B} \frac{1-r_k}{1-r_kx^{-1}} \prod_{k=1}^C (1-\beta_k + \beta_k x^{-1}),
    }
    where
    \begin{itemize}
        \item $G$ is a generating function of an $\amcm$ sequence with $G(1) \leq 1$;
        \item $q_k, r_k \in (0,1), \beta_k \in (0,1]$;
        \item $A,B,C$ are nonnegative integers.
    \end{itemize}
    If $B \leq C$, then there is $y \in \N$ and a diagonal random walk $(X_n, Y_n)$ starting from $(0,y)$ such that
    \formula[eq:main:E]{
    F(x) = \E^{(0,y)} \bigl[ x^{(X(\tau_0) - y)/2} \ind_{\tau_0 < \infty} \bigr]
    }
    for $|x|=1$. One can take $y=2A+B+C+1$.
    
    Conversely, every nonzero generating function of the form \eqref{eq:main:E} admits a representation of the form \eqref{eq:main:F} with $B \leq C$.
\end{theorem}

In the proof of the first part of the theorem we construct an appropriate random walk. Following Bondesson's approach, we introduce a barrier at level $y-1$, preventing the walk from returning to level $y$ after its first visit to level $y-1$. This allows us to construct the upper and lower parts of the walk separately. On levels $n\geq y$, we realise the $\amcm$ factor by iterating a transformation of $\amcm$ generating functions. On levels $n < y$, we realise all geometric factors using two-level blocks separated by further barriers. Condition $B \leq C$ reflects the use of the two-point factors in this construction. The remaining $C-B$ two-point factors can be realised with simple one-level blocks.

Bondesson's theorem for one-dimensional random walks can be recovered from our result. If $F$ is given by \eqref{eq:main:F} with $G$ being the generating function a one sided $\amcm$ sequence (hence, a completely monotone sequence) and $B=C=0$, then the condition $B \leq C$ is satisfied automatically. Our construction gives $p_{n,2}=p_{n,3}=0$ for every $n \in \N$. Hence, one can find a random walk $(X_n, Y_n)$ which starts from some $(0,2A+1)$ and moves only to the right (with NE and SE jumps). Consequently, we have $X(\tau_0) = \tau_0$ and Theorem 2 in \cite{bondesson} is recovered.

In the proof of the second part of Theorem \ref{thm:main} (see \eqref{deg:bounds}), we actually prove sharper bounds on the number of geometric factors in Theorem 1.1 in \cite{fpl}. We state this result below.

\begin{theorem}[Extension of Theorem 1.1 in \cite{fpl}]
\label{thm:main:1.5}
    Let $(X_n, Y_n)$ be a diagonal starting from $(0,y)$, where $y \geq 1$, such that $\pr(\tau_0<\infty)>0$. The distribution of $(X(\tau_0)-y)/2$ restricted to $\{\tau_0<\infty\}$ is given by a convolution of
    \begin{itemize}
        \item $A$ geometric sequences,
        \item $B$ reflected geometric sequences,
        \item at most $y$ reflected two-point sequences,
        \item an $\amcm$ sequence,
    \end{itemize}
    where $K_y = \lfloor (y-1)/2 \rfloor$ and $A,B$ are nonnegative integers such that
    \formula{
    &
    A \leq \min\biggl\{ K_y,  \sum_{k=2}^{y-1} \ind_{p_{k,4}>0} \biggr\}, &
    B \leq \min\biggl\{ K_y, \sum_{k=2}^{y-1} \ind_{p_{k,3}>0} \biggr\}.
    &
    }
\end{theorem}

\begin{remark}
Theorems~\ref{thm:main} and~\ref{thm:main:1.5} have a symmetric counterpart. If we consider the reflected random walk $\smash{(\tilde X_n, \tilde Y_n) = (-X_n,Y_n)}$ and normalise $\smash{(\tilde X(\tau_0)+y)/2}$, then the corresponding characteristic function is $\smash{\tilde F(x) = F(x^{-1})}$. Since $G(x^{-1})$ is again a generating function of an $\amcm$ sequence, \eqref{eq:main:F} reads
\formula{
\tilde F(x) = \tilde G(x) \prod_{k=1}^A \frac{1-q_k}{1-q_k x^{-1}} \prod_{k=1}^B \frac{1-r_k}{1-r_k x} \prod_{k=1}^C (1-\beta_k +\beta_k x).
}
With the notation above, the sufficient condition remains $B\leq C$, and one can still take $y=2A+B+C+1$. After relabelling the numbers of geometric and reflected geometric factors as $A$ and $B$, respectively, these become $A\leq C$ and $y=A+2B+C+1$. In the reflected version of Theorem~\ref{thm:main:1.5}, the reflected two-point sequences are replaced by two-point sequences supported on $\{0,1\}$. The bounds on the relabelled numbers $A$ and $B$ remain in their original form when expressed in terms of the reflected transition probabilities.
\end{remark}

\begin{remark}
By translation, Theorems~\ref{thm:main} and~\ref{thm:main:1.5} have analogous formulations for walks starting from arbitrary $(x,y) \in \Z \times \N$ and stopped upon hitting any integer level $a<y$. The corresponding formulas use $X_n-x$ and $y-a$ in place of $X_n$ and $y$, respectively.
\end{remark}

Clearly, the factorisation of the generating function of the first passage location is not unique. For example, consider a walk starting from $(0,1)$ that jumps only to the right. At level 1, it jumps NE with probability $q$ and SE with probability $1-q$ for some $q \in (0,1)$. At level 2, it jumps SE with probability one. The corresponding generating function is $(1-q)/(1-qx)$. Clearly, $B=C=0$ but since geometric sequences are $\amcm$, one can put
\formula{
&\text{either}\quad  G(x) = \frac{1-q}{1-qx},\, A=0 & \text{or}\quad G(x) = 1, \, A=1, \, q_1=q. &
}

A natural approach to this problem is to study representations that are maximal with respect to absorption into the $\amcm$ factor. We say that the representation \eqref{eq:main:F} of $F$ is \textit{maximal} if the functions
\formula{
& G(x)\frac{1-q_k}{1-q_kx}, & G(x) \frac{1-r_k}{1-r_kx^{-1}} &
}
are not generating functions of $\amcm$ sequences for any $1 \leq k \leq A$ and $1 \leq k \leq B$, respectively. Our next result concerns the maximal representations.

\begin{theorem}
\label{thm:main:2}
    If $F$ is a nonzero generating function of a first passage location \eqref{eq:main:E}, then every maximal representation of $F$ satisfies $B \leq C$.
\end{theorem}

Although the factorisation \eqref{eq:main:F} need not be unique, realisability of a function $F$ can be determined from any maximal representation. In particular, if $B>C$ in its maximal representation, then $F$ cannot be realised as a generating function of first passage location for any starting level~$y$.

\subsection{Organisation}
The paper is organised as follows. In Section~\ref{sec:rw}, we define our random walk, the corresponding generating functions, and recall their recurrence relations. In Section~\ref{sec:amcm}, we study $\amcm$ sequences and their generating functions. Section~\ref{sec:pick} is devoted for discussion on Pick functions. In Section~\ref{sec:tools}, we develop the main tools used in this paper. Section~\ref{sec:proofs} is devoted to proving the main results. Finally, Section~\ref{sec:discussion} presents examples and discusses possible extensions of our results to other classes of random walks.

\section{Random walks}
\label{sec:rw}

We consider a time-homogeneous Markov chain $(X_n,Y_n)$ on $\Z^2$ with diagonal jumps $(\pm1,\pm1)$. Its transition probabilities depend only on the current level $Y_n$. Following \cite{fpl}, we number the jump directions according
to the quadrants of the coordinate plane:
\formula{
p_{y,1} = \pr(X_{n+1} = X_n+1, \, Y_{n+1}=y+1 \, | \, X_n =x, Y_n =y ), \\
p_{y,2} = \pr(X_{n+1} = X_n-1, \, Y_{n+1}=y+1 \, | \, X_n =x, Y_n =y ), \\
p_{y,3} = \pr(X_{n+1} = X_n-1, \, Y_{n+1}=y-1 \, | \, X_n =x, Y_n =y ), \\
p_{y,4} = \pr(X_{n+1} = X_n+1, \, Y_{n+1}=y-1 \, | \, X_n =x, Y_n =y ), 
}
where $\sum_{k=1}^4 p_{y,k}=1$ for every $y \in \Z$. We write $\smash{\E^{(x,y)}}$ for expectation when the walk starts from $(x,y)$. By invariance under horizontal translations, it suffices to consider starting points of the form $(0,y)$.

Let $a \in \Z$. Our main interest is the \textit{first passage location}. It is defined as the position of the first coordinate of the walk while the second coordinate first reaches level $a$, that is, $X(\tau_a)$, where $\tau_a = \min\{n \geq 0: Y_n =a\}$ is the first passage time. We consider the rescaled passage location $(X(\tau_a) - y+a)/2$, So that it takes integer values. We are interested in passage to level $a=0$ from a starting level $y \in \N$.

We define the generating function of the rescaled first passage location as
\formula{
F_y^{a,b}(x) = \E^{(0,y)} \bigl[ x^{\frac{1}{2}(X(\tau_a) - y+a)} \ind_{\tau_a < \tau_b} \bigr]
}
for $|x|=1$ and $a,b,y \in \Z$ such that $\min\{a,b\} \leq y \leq \max\{a,b\}$. This corresponds to the distributions of the locations of the random walk at the first time of hitting level $a$ from $(0,y)$, restricted to the event that the walk reaches level $a$ before level $b$. For a random walk without a barrier, we define
\formula[eq:F:limit]{
F_y^{a,\infty}(x) = \lim_{b \to \infty} F_y^{a,b}(x) = \E^{(0,y)} \bigl[ x^{\frac{1}{2}(X(\tau_a) - y+a)} \ind_{\tau_a < \infty} \bigr].
}
In particular, the last equality follows from the dominated convergence theorem. As in Eq. (2.13) in \cite{fpl}, for $y \in \N$ and $a=0$, we have
\formula[eq:F:infty]{
F_y^{0,\infty}(x) = \frac{(p_{y,3}x^{-1}+p_{y,4}) F_{y-1}^{0,y}(x)}{1-(p_{y,1}x+p_{y,2})F_{y+1}^{y,\infty}(x) - (p_{y,3}x^{-1}+p_{y,4})F_{y-1}^{y,0}(x)}.
}
By Eqs. (2.4) and (2.12) in \cite{fpl}, for $y< b$,
\formula[eq:F:b:1]{
& F_y^{y-1,b}(x) = \frac{p_{y,3}x^{-1}+p_{y,4}}{1 - (p_{y,1}x+p_{y,2})F_{y+1}^{y,b}(x)}, & F_b^{b-1,b}(x) = 0; & 
}
The quotient identities above are understood wherever their denominators are nonzero. Observe that an analogous recurrence holds if we replace $b$ with $\infty$:
\formula[eq:F:b]{
F_y^{y-1,\infty}(x) = \frac{p_{y,3}x^{-1}+p_{y,4}}{1 - (p_{y,1}x+p_{y,2})F_{y+1}^{y,\infty}(x)}.
}
In this case, however, there is no terminal condition. Finally, observe that if $p_{y-1,1}=p_{y-1,2} = 0$, then $\smash{F_{y-1}^{y,0}} \equiv 0$ and $\smash{F_y^{0,\infty}}$ admits the following factorisation:
\formula[eq:F:factors]{
F_y^{0,\infty}(x) = F_y^{y-1, \infty}(x) \cdot F_{y-1}^{0,y}(x).
}

For a detailed discussion of generating functions of first passage locations, we refer the reader to Section 2.2 in \cite{fpl}.

\section{\texorpdfstring{$\amcm$}{AMCM} sequences and their generating functions}
\label{sec:amcm}

Recall that a (one-sided) sequence $(a(k): k \geq 0)$ is \textit{completely monotone} if $(-1)^n \Delta^n a(k) \geq 0$ for every $n,k \geq 0$, where $\Delta a(k) = a(k+1)-a(k)$ and $\Delta^n a(k) = \Delta \Delta^{n-1} a(k)$. The class of $\amcm$ sequences can be considered as a two-sided extension of completely monotone sequences. We say that a (two-sided) sequence $(a(k): k \in \Z)$ is $\amcm$ (absolutely monotone then completely monotone) if the sequences $(a(k): k \geq 0)$ and $(a(-k):k \geq 0)$ are completely monotone. We are mainly interested in probability and subprobability sequences. Therefore, we additionally assume that the mass of an $\amcm$ sequence $(a(k))$ is given by
\formula{
\sum_{k \in \Z} a(k) = M \in (0,1].
}
In particular, $a(k) \to 0$ as $k \to \pm\infty$.

The class of $\amcm$ sequences originated in \cite{bss-2}, where it played a key role in the characterisation of bell-shaped sequences. Its continuous analogue, the class of $\amcm$ functions, was also used to characterise bell-shaped functions in \cite{bsf-1, bsf-2}. In recent works of the author \cite{toeplitz-1, toeplitz-2}, a generalisation of this class was involved in the study of spectral theory of Toeplitz matrices and operators.

We focus on describing generating functions of $\amcm$ sequences. Hausdorff's classical theorem states that completely monotone sequences can be characterised as moment sequences. This property gives rise to the following characterisation theorem for $\amcm$ sequences.

\begin{theorem}
\label{thm:G}
Let $G$ be a holomorphic function on $\C\setminus[0,\infty)$ such that
\formula[eq:G:exp:lim]{
\lim_{t\to0}G(e^{it})=M\in(0,1].
}
We set $G(1)=M$. The following conditions are equivalent:
\begin{enumerate}[label = (\alph*)]
    \item\label{it:a:G} There exists an $\amcm$ sequence $(a(k))$ with $\sum_{k\in\Z}a(k)=M$ such that, for $|x|=1$,
    \formula{
    G(x)=\sum_{k\in\Z}a(k)x^k.
    }

    \item \label{it:b:G} There exist finite nonnegative Borel measures $\mu_\pm$ on $[0,1)$ with the same total mass
    \formula{
    m=\mu_+[0,1)=\mu_-[0,1)
    }
    such that, for $x\in\C\setminus[0,\infty)$,
    \formula[eq:G:int:1]{
    G(x)=\int_{[0,1)}\frac{1}{1-sx}\,\mu_+(ds)+\int_{[0,1)}\frac{1}{1-sx^{-1}}\,\mu_-(ds)-m,
    }
    or, equivalently,
    \formula[eq:G:int:2]{
    G(x)=\int_{(0,1)}\frac{sx}{1-sx}\,\mu_+(ds)+\int_{(0,1)}\frac{x}{x-s}\,\mu_-(ds)+\mu_-\{0\}.
    }

    \item\label{it:c:G} There exist a constant $c\in\R$ and a Borel function $\ph$ taking values in $[-1,0]$ on $(0,1)$ and in $[0,1]$ on $(1,\infty)$ such that, for $x\in\C\setminus[0,\infty)$,
    \formula[eq:G:exp]{
    G(x)=\exp\biggl(c+\int_0^\infty\biggl(\frac{1}{s-x}-\frac{s}{s^2+1}\biggr)\ph(s)\,ds\biggr).
    }
\end{enumerate}
The sequence in \ref{it:a:G} is uniquely determined by $G$. The measures in \ref{it:b:G} satisfy, for $k\geq0$,
\formula{
a(k)=\int_{[0,1)}s^k\,\mu_+(ds),
\qquad
a(-k)=\int_{[0,1)}s^k\,\mu_-(ds).
}
In particular, $m=a(0)\in(0,M]$. Moreover, for almost every $s>0$,
\formula[eq:pick:phi]{
\ph(s)=\frac{1}{\pi}\lim_{t\to0^+}\operatorname{Arg}G(s+it).
}
\end{theorem}

\begin{proof}
The implications \ref{it:a:G}$\Rightarrow$\ref{it:b:G} and \ref{it:a:G}$\Rightarrow$\ref{it:c:G} follow from Section~2.3 of \cite{bss-2} (see Lemma~2.2 therein). Conversely, the same results show that either \ref{it:b:G} or \ref{it:c:G} implies that $G$ is the generating function of a nonzero $\amcm$ sequence $(a(k))$ converging to zero as $k\to\pm\infty$. Indeed, the boundary condition required in Lemma~2.2 follows from our assumption:
\formula{
\lim_{t\to0}(e^{it}-1)G(e^{it})=0.
}
The moment formulas and \eqref{eq:pick:phi} also follow from these results.

It remains to verify that the resulting sequence is summable with mass $M$. Taking real parts in \eqref{eq:G:int:1}, we obtain, for $0<t<\pi$,
\formula{
    \re G(e^{it}) &= \int_{[0,1)} \frac{1-s \re e^{-it}}{|1-se^{it}|^2} \, \mu_+(ds) + \int_{[0,1)} \frac{1-s \re e^{it}}{|1-se^{-it}|^2} \, \mu_-(ds) - m \\
    &= \int_{[0,1)} \frac{1-s \re e^{it}}{|1-se^{it}|^2} \, (\mu_++\mu_-)(ds) - m \\
    &= \int_{[0,1)} \frac{1-s \cos t}{1-2s\cos t+s^2} \, (\mu_++\mu_-)(ds) - m \\ 
    &= \frac{1}{2} \int_{[0,1)} \biggl( 1+\frac{1-s^2}{1-2s\cos t+s^2} \biggr) \, (\mu_++\mu_-)(ds) - m.
}
Using $\mu_+[0,1)=\mu_-[0,1)=m$, we find that
\formula{
\operatorname{Re}G(e^{it})=\frac12\int_{[0,1)}\frac{1-s^2}{(1-s)^2+2s(1-\cos t)}\,(\mu_++\mu_-)(ds).
}
The integrand is nonnegative and increases to $(1+s)/(1-s)$ as $t\to 0^+$. On one hand, by the monotone convergence theorem,
\formula{
M = \lim_{t \to 0} G(e^{it}) = \frac{1}{2} \int_{[0,1)} \frac{1+s}{1-s} \, (\mu_++\mu_-)(ds).
}
On the other hand, by Tonelli theorem,
\formula{
\sum_{k\in\Z}a(k) &= \sum_{k =1}^\infty \int_{[0,1)} s^k \, (\mu_++\mu_-)(ds)+m \\
&= \int_{[0,1)} \frac{s}{1-s} \, (\mu_++\mu_-)(ds)+m \\
&= \frac{1}{2} \int_{[0,1)} \biggl( \frac{1+s}{1-s}-1\biggr) \, (\mu_++\mu_-)(ds)+m \\
&= \frac{1}{2} \int_{[0,1)} \frac{1+s}{1-s} \, (\mu_++\mu_-)(ds).
}
Therefore, 
\formula{
M = \sum_{k \in \Z} a(k).
}
Thus, \ref{it:a:G} holds, including the identity at $x=1$. The sequence is uniquely determined by the Fourier coefficients of $G$. Finally, $m\leq M$ since $m=a(0)\leq\sum_{k\in\Z}a(k)=M$,
\end{proof}

By $\mathscr{G}$ we denote the class of nonzero functions described in Theorem~\ref{thm:G}. Note that $\mathscr{G}$ is closed under inversion of arguments. More precisely, if $G \in \mathscr{G}$, then $\tilde G(x) = G(x^{-1})$ is also in $\mathscr{G}$. We now prove some basic properties of this class.

\begin{lemma}
\label{lem:G:monotone}
    If $G \in \mathscr{G}$, then exactly one of the following is true:
    \begin{enumerate}
        \item $G$ is constant on $(-\infty,0)$,
        \item $G$ is strictly monotone on $(-\infty,0)$,
        \item $G$ is strictly increasing on $(-\infty, -\alpha]$ and strictly decreasing on $[-\alpha,0)$, and it attains maximum at $-\alpha$ for some $\alpha>0$.
    \end{enumerate}
\end{lemma}

\begin{proof}
    Let $x \in (0,1)$ and denote $\tilde G(x) = G(x/(x-1))$. Using \eqref{eq:G:int:1}, we get
    \formula{
    \tilde G(x) = \int_{[0,1)} \frac{1-x}{1-x+sx} \, \mu_+(ds) + \int_{[0,1)} \frac{x}{x+s-sx} \, \mu_-(ds) - \mu_+[0,1).
    }

    If $\mu_+(0,1)=\mu_-(0,1)=0$, then both measures are supported at zero. Since they have equal total masses, $\tilde G$ is identically equal to $\mu_+\{0\}$ on $(0,1)$. Consequently, $G$ is constant on $(-\infty,0)$.
    
    Otherwise, at least one of $\mu_+(0,1)$ and $\mu_-(0,1)$ is positive. For each fixed $s \in (0,1)$, the functions
    \formula{
    & f_s(x)=\frac{1-x}{1-x+sx}, & g_s(x)=\frac{x}{x+s-sx} &
    }
    are strictly concave on $(0,1)$. Indeed,
    \formula{
    & f_s''(x)=-\frac{2s(1-s)}{(1-x+sx)^3}, & g_s''(x)=-\frac{2s(1-s)}{(x+s-sx)^3}. &
    }
    Both functions take values in $[0,1]$, so they are integrable with respect to the finite measures in the representation above. Integrating their concavity inequalities, and using the fact that at least one measure assigns positive mass to $(0,1)$, shows that $\tilde G$ is strictly concave.

    A strictly concave function on an interval is either strictly monotone or strictly increasing and then strictly decreasing. Finally, the map $x\mapsto x/(x-1)$ is a strictly decreasing bijection from $(0,1)$ onto $(-\infty,0)$. Thus, strict monotonicity of $\tilde G$ implies strict monotonicity of $G$. If $\tilde G$ has its maximum at $x_0$, then $G$ is strictly increasing up to $w_0=x_0/(x_0-1)$ and strictly decreasing thereafter.
\end{proof}

\begin{lemma}
\label{lem:G:limits}
    If $G \in \mathscr{G}$, then
    \formula{
    & \lim_{x \to 0^-} G(x) = \mu_-\{0\}, & \lim_{x \to -\infty} G(x) = \mu_+\{0\}. &
    }
\end{lemma}

\begin{proof}
    By \eqref{eq:G:int:2},
    \formula{
    G(x) = \int_{(0,1)} \frac{sx}{1-sx} \, \mu_+(ds) + \int_{(0,1)} \frac{x}{x-s} \, \mu_-(ds) + \mu_-\{0\}.
    }
    The limit at $0^-$ is a consequence of the dominated convergence theorem. Indeed, for $x \in (-1,0)$ and $s \in (0,1)$, we have
    \formula{
    & \bigg\lv \frac{sx}{1-sx} \bigg\rv < 1, & \bigg\lv \frac{x}{x-s} \bigg\rv < 1.
    }
    The limit at $-\infty$ follows by applying the first assertion to $\tilde G(x)=G(x^{-1})$, which belongs to $\mathscr{G}$ and has representing measures $\tilde\mu_+=\mu_-$ and $\tilde\mu_-=\mu_+$.
\end{proof}

\begin{lemma}
\label{lem:G:bounds}
    If $G \in \mathscr{G}$ is nonconstant, then for $x \in (-\infty,0)$,
    \formula{
    \min\{\mu_+\{0\}, \mu_-\{0\}\} < G(x) < a(0).
    }
    Moreover,
    \formula{
    0 <\sup_{x<0}G(x) < M.
    }
\end{lemma}

\begin{proof}
    The lower bound on $G$ follows from Lemmas \ref{lem:G:monotone} and \ref{lem:G:limits}. Let $t >0$. By \eqref{eq:G:int:1},
    \formula{
    G(-t) &= \int_{[0,1)} \frac{1}{1+st} \, \mu_+(ds) + \int_{[0,1)} \frac{1}{1+st^{-1}} \, \mu_-(ds) - a(0) \\
    &< a(0) + a(0)  - a(0) = a(0).
    }
    Here we used $\mu_+[0,1) = \mu_-[0,1)= a(0)$. This completes the proof of the first part.
    
    If $G$ is nonconstant, then $a(0) < M$. Otherwise, we would have $a(0) = \sum_{k \in\Z} a(k)$ and thus, $a(k)=0$ for $k \neq 0$. Moreover, using the integral representation \eqref{eq:G:int:1}, we deduce that
    \formula{
    G(-1) &= \int_{[0,1)} \frac{1}{1+s} \, \mu_+(ds) + \int_{[0,1)} \frac{1}{1+s} \, \mu_-(ds) - a(0)\\
    &> \frac{1}{2} a(0) + \frac{1}{2} a(0) - a(0) = 0.
    }
    It follows that
    \formula{
    0 < G(-1) \leq \sup_{x<0}G(x) \leq a(0) < M,
    }
    as claimed.
\end{proof}

\begin{lemma}
\label{lem:G:limits:2}
    If $G \in \mathscr{G}$, then the limits
    \formula{
    & \lim_{x \to 0^-} \frac{G(x)}{x}, & \lim_{x \to -\infty} xG(x) &
    }
    exist and belong to $[-\infty,0)$.
\end{lemma}

\begin{proof}
    If $\mu_-\{0\}>0$, then applying Lemma \ref{lem:G:limits} yields
    \formula{
    \lim_{x \to 0^-} \frac{G(x)}{x} = -\infty.
    }
    Assume that $\mu_-\{0\}=0$ and let $x \in (-1,0)$. Using \eqref{eq:G:int:2}, we find that
    \formula{
    \frac{G(x)}{x} = \int_{(0,1)} \frac{s}{1-sx} \, \mu_+(ds) + \int_{(0,1)} \frac{1}{x-s} \, \mu_-(ds).
    }
    Note that $|s/(1-sx)| \leq s$ for $s \in (0,1)$, $x \in (-1,0)$, and $\mu_+$ has finite moments of nonnegative orders. Applying the dominated convergence theorem to the first integral yields
    \formula{
    \lim_{x \to 0^-} \int_{(0,1)} \frac{s}{1-sx }\,\mu_+(ds) = \int_{(0,1)} s\,\mu_+(ds) < \infty.
    }
    To calculate the second integral, we use the monotone convergence theorem. Since $x \mapsto (s-x)^{-1}$ is a nonnegative strictly increasing function as $x \to 0^-$, we get
    \formula{
    \lim_{x\to 0^-} \int_{(0,1)} \frac{1}{x-s} \, \mu_-(ds) = -\int_{(0,1)} s^{-1} \, \mu_-(ds),
    }
    where the right-hand side may equal $-\infty$. Consequently, using the assumption that $\mu_-\{0\}=0$, we estimate
    \formula{
    \lim_{x \to 0^-} \frac{G(x)}{x} &= \int_{(0,1)} s \, \mu_+(ds) - \int_{(0,1)} s^{-1} \, \mu_-(ds) < \mu_+[0,1) - \mu_-[0,1)~=~0.
    }
    Finally, applying the first part to $\tilde G(x)=G(x^{-1})$ gives the second part of the assertion.
\end{proof}

\begin{lemma}
\label{lem:reG}
    If $G \in \G$, then $\re G(x) >0$ for $|x|=1$.
\end{lemma}

\begin{proof}
    Let $G \in \G$. By \eqref{eq:G:int:1}, we have
    \formula{
    G(x) = \int_{[0,1)}  \frac{1-s\overline{x}}{|1-sx|^2} \, \mu_+(ds) + \int_{[0,1)} \frac{1-s\overline{x}^{-1}}{|1-sx^{-1}|^2} \mu_-(ds) - a(0).    
    }
    If $|x|=1$, then $\overline x = x^{-1}$ and $|1-sx^{-1}| = |1-sx|$. Consequently,
    \formula{
    \re G(x) = \int_{[0,1)} \frac{1-s \re x}{|1-sx|^2} (\mu_++\mu_-)(ds) - a(0).
    }
    Since $2(1-s\re x) =|1-sx|^2 + 1-s^2$, we get
    \formula{
    \re G(x) &= \frac{1}{2} \int_{[0,1)} \biggl( 1 + \frac{1-s^2}{|1-sx|^2}\biggr)  (\mu_++\mu_-)(ds) - a(0) \\
    &= \frac{1}{2} \int_{[0,1)} \frac{1-s^2}{|1-sx|^2} (\mu_++\mu_-)(ds) + \frac{1}{2}(a(0) + a(0)) - a(0) \\
    &= \frac{1}{2} \int_{[0,1)} \frac{1-s^2}{|1-sx|^2} (\mu_++\mu_-)(ds)>0.
    }
    Above we used the property $\mu_+[0,1)=\mu_-[0,1)=a(0)$.
\end{proof}

\section{Pick functions}
\label{sec:pick}

A \textit{Pick function} is a function $P$ holomorphic in the upper half-plane $\{x \in \C: \im x > 0\}$ with nonnegative imaginary part. In other words, $\im f(x) \geq 0$ for $\im x >0$. Every Pick function admits the following integral representation
\formula[eq:pick:int:1]{
P(x) = Ax + B + \int_{\R} \biggl( \frac{1}{s-x}  - \frac{s}{s^2+1} \biggr) \, \mu(ds),
}
where $A \geq 0$, $B \in \R$ and $\mu$ is a nonnegative measure satisfying $\int_{\R} (s^2+1)^{-1} \, \mu(ds)<\infty$. Another classical result states that if $P$ is a Pick function not identically equal to zero, then $\log P$ is again a Pick function. It follows that $P$ has the following exponential representation
\formula[eq:pick:exp]{
P(x) = \exp\biggl( c + \int_{-\infty}^\infty \biggl( \frac{1}{s-x} - \frac{s}{s^2+1} \biggr) \, \ph(s) \, ds \biggr)
}
for $\im x>0$, where $c \in \R$ and $\ph$ is a Borel function that takes values in $[0,1]$. In fact, for a nonzero Pick function, representations \eqref{eq:pick:int:1} and \eqref{eq:pick:exp} are equivalent. For more details and discussion, see Chapter II of \cite{donoghue} or Section 2.1 in \cite{bss-2}.

In this paper, we focus on Pick functions that, for some $\alpha>0$, admit a holomorphic extension to $\C\setminus(\{-\alpha\}\cup[0,\infty))$ which is real-valued on $(-\infty,0)\setminus\{-\alpha\}$. This is equivalent to requiring that the measure $\mu$ in \eqref{eq:pick:int:1} be supported on $\{-\alpha\}\cup[0,\infty)$. Note that
\formula{
\frac{1}{s-x} - \frac{s}{s^2+1} = \frac{1}{s^2+1} \, \frac{1+sx}{s-x}.
}
Putting $\sigma(ds) = (1+s^2)^{-1}\mu(ds)$, one can rewrite \eqref{eq:pick:int:1} as
\formula{
    P(x) &= Ax + B + \int_{\{-\alpha\} \cup [0,\infty)} \frac{1+sx}{s-x} \, \sigma(ds) \\
    &= Ax + B - \frac{\sigma\{0\}}{x} - \frac{1-\alpha x}{\alpha+x} \sigma\{-\alpha\} + \int_{(0,\infty)} \frac{1+sx}{s-x} \, \sigma(ds) \\
    &= \bigl(B + \alpha \sigma\{-\alpha\} \bigr) + Ax - \frac{\sigma\{0\}}{x}  - \frac{\sigma\{-\alpha\}(1+\alpha^2)}{\alpha+x} + \int_{(0,\infty)} \frac{1+sx}{s-x} \, \sigma(ds).
}
Hence, in that case, we have another integral representation, given by
\formula[eq:pick:int:2]{
P(x) = Ax + B - \frac{C}{x} - \frac{D}{\alpha+x} + \int_{(0,\infty)} \frac{1+sx}{s-x} \, \sigma(ds),
}
where $B \in \R$, $A,C,D \geq 0$ and $\sigma$ is a finite measure on $(0,\infty)$.

Pick functions of that form are closely related to generating functions of $\amcm$ sequences. More precisely, if $G \in \G$, then $(1-x^{-1})G(x)$, which is a generating function of $(-\Delta a(k))$, is a Pick function given by \eqref{eq:pick:int:2} with $A=D=0$ and $\sigma\{1\}=0$. For details, see Step 2 in the proof of Lemma 2.2 in \cite{bss-2}. 

Below we collect some auxiliary results for the above subclass \eqref{eq:pick:int:2} of Pick functions.

\begin{lemma}
\label{lem:pick:constants}
    If $P$ is a Pick function given by \eqref{eq:pick:int:2}, then
    \formula{
    & A=\lim_{x \to -\infty} \frac{P(x)}{x}, && C=-\lim_{x \to 0^-} xP(x), & D=-\lim_{x \to -\alpha} (x+\alpha)P(x). &
    }
\end{lemma}

\begin{proof}
    First assume that $x\in(-\infty,-1)$ and observe that for $s>0$ we have
    \formula{
    \bigg\lv \frac{1+sx}{x(s-x)} \bigg\rv \leq \frac{1+s|x|}{|x|(s+|x|)} < 1.
    }
    Since $\sigma$ is finite, the dominated convergence theorem therefore gives
    \formula{
    \lim_{x \to -\infty} \frac{P(x)}{x} = A + \int_{(0,\infty)} \lim_{x\to -\infty} \frac{1+sx}{x(s-x)} \, \sigma(ds) = A.
    }
    Similarly, for $x \in (-1,0)$ and $s>0$ we have
    \formula{
    \biggl\vert \frac{x(1+sx)}{s-x} \biggr\vert \leq \frac{|x|(1+s|x|)}{s+|x|} < 1.
    }
    Again by the dominated convergence theorem, we get
    \formula{
    \lim_{x\to0^-} xP(x) = -C + \int_{(0,\infty)} \lim_{x\to 0^-} \frac{x(1+sx)}{s-x} \, \sigma(ds) = -C.
    }
    Finally, let $|x+\alpha|<\alpha/2$. Then, for $s>0$,
    \formula{
    \bigg\lv \frac{(\alpha+x)(1+sx)}{s-x}\bigg\rv \leq \frac{|\alpha+x|(1+s|x|)}{s+|x|} < \frac{\alpha}{2} \, \frac{2+3\alpha s}{2s + \alpha} \leq \max\biggl\{ 1, \frac{3\alpha^2}{4} \biggr\}.
    }
    Consequently,
    \formula{
    \lim_{x \to -\alpha} (x+\alpha) P(x) = -D +\int_{(0,\infty)} \lim_{x\to -\alpha} \frac{(\alpha+x)(1+sx)}{s-x} \, \sigma(ds) = -D,
    }
    and the proof is complete.
\end{proof}

\begin{lemma}
\label{lem:pick:limits}
    Let $P$ be a Pick function given by \eqref{eq:pick:int:2}. If $C=0$, then
    \formula{
    & \lim_{x \to 0^-} P(x) = B - \frac{D}{\alpha} + \int_{(0,\infty)} s^{-1} \, \sigma(ds).
    }
    If $A=0$, then
    \formula{
    &\lim_{x \to -\infty} P(x) = B - \int_{(0,\infty)} s\,\sigma(ds).
    }
\end{lemma}

\begin{proof}
    If $C=0$, the integral representation of $P$ yields
    \formula{
    \lim_{x\to0^-} P(x) = B - \frac{D}{\alpha} + \lim_{x\to0^-} \int_{(0,\infty)} \frac{1+sx}{s-x} \, \sigma(ds).
    }
    Denote
    \formula{
    g_s(x) = \frac{1+sx}{s-x}
    }
    and observe that $g_s(-1)$ is integrable: for $s>0$, we have $|g_s(-1)| = |(1-s)/(1+s)| \leq 1$. Moreover, the function $x \mapsto g_s(x) - g_s(-1)$ is nonnegative for $x \in (-1,0)$ and increases as $x \to 0$. By the monotone convergence theorem,
    \formula{
    \lim_{x \to 0^-} \int_{(0,\infty)} g_s(x) \, \sigma(ds) &= \lim_{x \to 0^-}  \int_{(0,\infty)} (g_s(x)-g_s(-1)) \, \sigma(ds) + \int_{(0,\infty)} g_s(-1) \, \sigma(ds) \\
    &=\int_{(0,\infty)} (s^{-1} - g_s(-1)) \, \sigma(ds) + \int_{(0,\infty)} g_s(-1) \sigma(ds) \\ &= \int_{(0,\infty)} s^{-1} \, \sigma(ds).
    }
    If $A=0$, a similar argument is applied to $g_s(-1)-g_s(x)$ as $x \to -\infty$.
\end{proof}

\section{Tools}
\label{sec:tools}

\subsection{Transformation of generating functions}
In this section we develop the main tools for proving Theorem \ref{thm:main}. Before that, we note some useful identities for $\alpha> 0$:
\formula[]{
    \log x &= \int_{-\infty}^0 \biggl( \frac{1}{s-x} - \frac{s}{s^2+1} \biggr) \, ds, \label{eq:x} \\
    \log (1-x) &= -\int_1^\infty \biggl( \frac{1}{s-x} - \frac{s}{s^2+1} \biggr) \, ds + \log \sqrt{2}, \label{eq:1-x} \\
    \log(1-x^{-1}) &= \int_0^1 \biggl( \frac{1}{s-x} - \frac{s}{s^2+1} \biggr) \, ds + \log \sqrt{2}, \label{eq:1-1/x} \\
    \log (1+\alpha x^{-1}) &= -\int_{-\alpha}^{0} \biggl( \frac{1}{s-x} - \frac{s}{s^2+1} \biggr) \, ds + \log \sqrt{\alpha^2+1}, \label{eq:1+a/x} \\
    \log (1+\alpha^{-1} x) &= \int_{-\infty}^{-\alpha} \biggl( \frac{1}{s-x} - \frac{s}{s^2+1} \biggr) \, ds + \log \sqrt{\alpha^{-2}+1}.  \label{eq:1+x/b}
    }
These identities hold for $\im x>0$ with the principal branches of the logarithms, and extend analytically wherever the corresponding branches are defined.

\begin{proposition}
\label{prop:transformation}
    Let $G \in \G$ be nonconstant and denote $g = \sup_{x<0}G(x)$. There are $\tilde G \in \G$ and $\beta \in [0,1]$ such that
    \formula{
    1-\frac{g}{G(x)} = (1-g) T_\beta(x) \tilde{G}(x),
    }
    where
    \formula[eq:T]{
    T_\beta(x) = (\beta x^{-1} + 1-\beta)((1-\beta)x+\beta).
    }
\end{proposition}

\begin{proof}
    By Lemma \ref{lem:G:monotone}, $G$ is either strictly monotone on $(-\infty,0)$, or it attains a maximum. First, we assume that $G$ is increasing on $(-\infty,-\alpha]$ and decreasing on $[-\alpha,0)$ for some $\alpha \in (0,\infty)$. Then $g = G(-\alpha) \in (0,M)$ by Lemma \ref{lem:G:bounds}. Our argument is divided into four key steps.
    
    \textit{Step 1.} We define the following auxiliary function
    \formula[eq:P:1]{
    P(x) = \frac{1}{(1+\alpha x^{-1})(1-x)G(x)},
    }
    which is analytic on $\C \setminus(\{-\alpha\}\cup[0,\infty))$. Since $G \in \G$, it has an exponential representation \eqref{eq:G:exp}. Using \eqref{eq:1-x} and \eqref{eq:1+a/x}, we find that $P$ has the following exponential representation:
    \formula[eq:P:2]{
    P(x) = \exp\biggl( c + \int_{-\infty}^\infty \biggl( \frac{1}{s-x} - \frac{s}{s^2+1} \biggr) (\ind_{(-\alpha,0)}(s)+\ind_{(1,\infty)}(s)-\varphi(s)) \, ds \biggr),
    }
    where $\ph \equiv 0$ on $(-\infty,0)$, $\ph \in [-1,0]$ on $(0,1)$, $\ph \in [0,1]$ on $(1,\infty)$, and $c\in \R$. Consequently, the function $\ind_{(-\alpha,0)}(s)+\ind_{(1,\infty)}(s)-\varphi(s)$ takes values in $[0,1]$. Therefore, $P$ is a Pick function analytic on $\C \setminus (\{-\alpha\} \cup [0,\infty))$. This allows us to rewrite it in the integral representation \eqref{eq:pick:int:2}:
    \formula[eq:P:3]{
    P(x) = Ax + B - \frac{C}{x} - \frac{D}{x+\alpha} + \int_{(0,\infty)} \frac{1+sx}{s-x} \, \sigma(ds),
    }
    where $A,C,D \geq 0$, $B \in \R$, and $\sigma$ is a finite measure.

    \textit{Step 2.} We now investigate \eqref{eq:P:3}. By Lemmas \ref{lem:G:limits:2} and \ref{lem:pick:constants}, we have
    \formula{
    A&= \lim_{x \to -\infty} \frac{P(x)}{x} = \lim_{x \to -\infty} \frac{1}{(1+\alpha x^{-1})(1-x)xG(x)} = 0, \\
    C&= -\lim_{x \to 0^-} xP(x) = -\lim_{x \to 0^-} \frac{x}{(1+\alpha x^{-1})(1-x) G(x)} = 0.
    }
    Moreover, since $G(-\alpha)=g$, we have
    \formula{
    D = -\lim_{x \to -\alpha} (x+\alpha) P(x) = -\lim_{x\to-\alpha} \frac{x}{(1-x)G(x)} = \frac{\alpha}{g(1+\alpha)}.
    }
    Note that $P$ is positive on $(-\infty,-\alpha)$ and negative on $(-\alpha,0)$. It is a consequence of \eqref{eq:P:1} and Lemma \ref{lem:G:bounds}, which states that $G$ is positive on $(-\infty,0)$. Using this property and Lemma \ref{lem:pick:limits}, we get
    \formula[]{
    0 \leq \lim_{x \to -\infty} P(x) &= B - \int_{(0,\infty)} s \, \sigma(ds), \label{eq:B1} \\
    0 \geq \lim_{x \to 0^-} P(x) &= B - \frac{1}{g(1+\alpha)} + \int_{(0,
    \infty)} s^{-1} \, \sigma(ds). \label{eq:B2}
    }
    In particular, the above integrals are finite. Hence, we have proved that
    \formula[eq:P:4]{
    P(x) = B - \frac{\alpha}{g(1+\alpha)(x+\alpha)} + \int_{(0,\infty)} \frac{1+sx}{s-x} \, \sigma(ds),
    }
    where $B \geq 0$ and $\sigma$ has finite moments of orders $-1,0,1$.

    \textit{Step 3.} Let
    \formula{
    \tilde P (x) = \frac{1-x^{-1}}{1+\alpha x^{-1}} \biggl( 1 -\frac{g}{G(x)}\biggr) = \frac{x-1}{x+\alpha} - g(1-x)(1-x^{-1}) P(x).
    }
    Using \eqref{eq:P:4} and the identities
    \formula{
    -(1-x)(1-x^{-1}) &= x+x^{-1}-2, \\
    -(1-x)(1-x^{-1})\frac{1+sx}{s-x} &= \frac{1+sx}{s-x} \, \frac{(1-s)^2}{s} - sx + \frac{1}{sx}, \\
    (1-x)(1-x^{-1})\,\frac{\alpha}{(1+\alpha)(x+\alpha)} &= -\frac{\alpha}{1+\alpha} - \frac{1}{x} \, \frac{1}{1+\alpha} + \frac{1+\alpha}{x+\alpha}, 
    }
    we rewrite $\tilde P$ as
    \formula[eq:tilde:P]{
    \tilde P(x) = \tilde A x + \tilde B - \frac{\tilde C}{x} + \int_{(0,\infty)} \frac{1+sx}{s-x} \, \tilde \sigma(ds),
    }
    where $\tilde \sigma(ds) = g(1-s)^2 s^{-1} \sigma(ds)$ and
    \formula{
    & \tilde A = g\biggl(B- \int_{(0,\infty)} s \, \sigma(ds)\biggr),
    & \tilde B =  1-2Bg - \frac{\alpha}{1+\alpha}, &&
    \tilde C = \frac{1}{1+\alpha} -g \biggl(B+\int_{(0,\infty)} s^{-1} \, \sigma(ds) \biggr), &
    }
    Note that since $\sigma$ has finite moments of orders $-1,0,1$, the constants $\tilde A, \tilde C$ are well defined and $\tilde \sigma$ is a finite measure on $(0,\infty)$. Moreover, by \eqref{eq:B1} and \eqref{eq:B2}, we find that $\tilde A, \tilde C \geq 0$. Consequently, $\tilde P$ defines a Pick function which is analytic on $\C \setminus[0,\infty)$.

    \textit{Step 4.} Observe that if $\beta = \alpha/(1+\alpha)$, then $\beta \in (0,1)$ and
    \formula{
    T_\beta(x) = \frac{(x+\alpha)^2}{x(1+\alpha)^2}.
    }
    We denote
    \formula{
    \tilde G(x) = \frac{1-g/G(x)}{(1-g)T_\beta(x)} = \frac{(1+\alpha)^2}{\alpha(1-g)} \, \frac{\tilde P(x)}{(1-x^{-1})(1+\alpha^{-1} x)}.
    }
    Our goal is to show that $\tilde G \in \G$. The exponential representation of the Pick function $\tilde P$ reads
    \formula[eq:P:exp:case:1:2]{
        \tilde P(x) = \exp\biggl( \tilde c + \int_{-\infty}^\infty \biggl( \frac{1}{s-x} - \frac{s}{s^2+1} \biggr) \, \tilde \ph(s) \, ds \biggr),
    }
    where $\tilde c \in \R$ and $\tilde \ph$ is given by \eqref{eq:pick:phi}. Since $G(x) < g$ for $x \in (-\infty,-\alpha)\cup(-\alpha,0)$, we have $\tilde P(x)<0$ for $x \in (-\infty,-\alpha)$ and $\tilde P(x) > 0$ for $x \in (-\alpha,0)$. It follows that $\tilde \ph \equiv 1$ on $(-\infty,-\alpha)$ and $\tilde \ph \equiv0$ on $(-\alpha,0)$. Hence, we write $\tilde \ph(s) = \ind_{(-\infty,-\alpha)}(s) + \psi(s)$, where $\psi \in [0,1]$ on $(0,\infty)$ and $\psi \equiv 0$ on $(-\infty,0)$. Using \eqref{eq:1-1/x} and \eqref{eq:1+x/b}, we obtain
    \formula{
    \tilde{G}(x) = \frac{(1+\alpha)^2}{\alpha(1-g)} \, \exp\biggl( d + \int_{-\infty}^\infty \biggl( \frac{1}{s-x} - \frac{s}{s^2+1} \biggr) \, (\psi(s) - \ind_{(0,1)}(s)) \, ds \biggr)
    }
    for some $d \in \R$. By the continuity of $G$ and $T_\beta$ at $1$, and since $G(1)=M$ and $T_\beta(1)=1$, we have
    \formula{
    \lim_{t\to0}\tilde G(e^{it})=\frac{1-g/M}{1-g}\in(0,1],
    }
    where we used $0<g<M\leq1$. Together with the exponential representation obtained above, Theorem~\ref{thm:G} implies that $\tilde G\in\G$. This completes the proof in the case where $G$ attains its maximum on $(-\infty,0)$.

    Now assume that $G$ is strictly increasing on $(-\infty,0)$. Hence, $g = \lim_{x\to0^-}G(x)$ and we put $\alpha=0$ in the proof for the previous case. Then \eqref{eq:P:1} reads
    \formula{
    P(x) = \frac{1}{(1-x)G(x)},
    }
    so $P$ is a Pick function analytic on $\C \setminus [0,\infty)$. In its integral representation \eqref{eq:P:3}, we have $A=C=D=0$. The integral representation \eqref{eq:tilde:P} can be rewritten by putting $\alpha =0$ in the formulas for $\tilde A$ and $\tilde B$. Moreover, by Lemma \ref{lem:pick:limits}, we have
    \formula{
    \tilde C = 1 - g \biggl(B+\int_{(0,\infty)} s^{-1}\, \sigma(ds) \biggr) = 1- g \lim_{x \to 0^-} P(x) = 1-g \cdot \frac{1}{g}=0.
    }
    Therefore $\tilde P$ is a Pick function analytic on $\C \setminus [0,\infty)$. Since $T_0(x) = x$, in step 4, we define
    \formula{
    \tilde G(x) = \frac{1-g/G(x)}{x(1-g)} = \frac{1}{1-g} \, \frac{\tilde P(x)}{x(1-x^{-1})}.
    }
    Note that $\tilde P(x)<0$ for all $x <0$. It follows that in the exponential representation of $\tilde P$ we have $\tilde \ph(s) = \ind_{(-\infty,0)}(s) + \psi(s)$, where $\psi \in [0,1]$ on $(0,\infty)$. Consequently, by \eqref{eq:x} and \eqref{eq:1-1/x}, we arrive at
    \formula{
    \tilde G(x) &= \frac{1}{1-g} \exp\biggl( d + \int_{-\infty}^\infty \biggl( \frac{1}{s-x} - \frac{s}{s^2+1} \biggr) \, (\psi(s) - \ind_{(0,1)}(s)) \, ds \biggr)
    }
    for some $d \in \R$. Moreover, $\tilde G(1) = (1-g/M)/(1-g) \in (0,1]$. Thus, $\tilde G \in \G$.

    The case when $G$ is strictly decreasing on $(-\infty,0)$ follows from applying the assertion to $H(x) = G(x^{-1})$ which is strictly increasing on $(-\infty,0)$.
\end{proof}


\begin{corollary}
\label{cor:recurrence}
    Let $G_y \in \G$ and denote $g_y = \sup_{x<0}G_y(x)$. There are constants $\beta_n \in [0,1]$ and functions $G_n \in \G \cup \{0\}$ with $g_n = \sup_{x<0}G_n(x)$, such that for every $n \geq y$, we have
    \formula[eq:Gn:recurrence]{
    G_n(x) = \frac{g_n}{1-(1-g_n)T_{\beta_{n+1}}(x) G_{n+1}(x)},
    }
    where $T_\beta$ is given by \eqref{eq:T}.
\end{corollary}

\begin{proof}
    If $G_y$ is nonconstant, by Proposition~\ref{prop:transformation}, there are $\beta_{y+1} \in [0,1]$ and $G_{y+1} \in \G$ such that
    \formula[eq:from:prop]{
    1-\frac{g_y}{G_y(x)} = (1-g_y) T_{\beta_{y+1}}(x) G_{y+1}(x).
    }
    Hence,
    \formula{
    G_{y}(x) = \frac{g_y}{1-(1-g_y) T_{\beta_{y+1}}(x) G_{y+1}(x)}.
    }
    Repeating this argument whenever $G_n$ is nonconstant yields either an infinite sequence or a first index $N\geq y$ for which $G_N$ is constant. In the latter case, including $N=y$, we set $G_n\equiv0$ and $\beta_n=0$ for all $n> N$. The recurrence then holds also for every $n\geq N$.
\end{proof}

\subsection{Truncated recursions and convergence}
Let $(G_n: n \geq y) \subset \G$ be the family of nonconstant functions from Corollary~\ref{cor:recurrence}. By $\beta_n \in [0,1]$ and $g_n \in (0,1)$ we denote the corresponding constants and suprema, respectively. For $N \geq y$, we consider the families of functions $(G_n^{N}:y\leq n \leq N)$ defined by the following recurrence relation for $y \leq n < N$,
\formula[eq:GN:recurrence]{
    & G_n^{N}(x) = \frac{g_n}{1-(1-g_n)T_{\beta_{n+1}}(x) G_{n+1}^{N}(x)}, & G_N^{N}(x) = g_N. &
}
Below we prove some properties of this family.

\begin{lemma}
\label{lem:GN:G}
    For $y \leq n \leq N$, we have $G_n^N \in \G$ and $0 < G_n^N(x) < 1$ for $x \in (-\infty,0)$.
\end{lemma}

\begin{proof}
    We use the class $\mathscr{Q}$ of rational functions with interlacing zeros and poles, introduced in \cite{fpl} (see Definition 4.4 therein). By Lemma 4.8 in \cite{fpl}, if $G \in \mathscr{Q}$ and $G(1) \leq 1$, then $G \in \G$. Lemma 4.6 in \cite{fpl} states that for $G \in \mathscr{Q}$ and $a,b,c \geq 0$ such that $a+b+c>0$ and $(a+b+c)G(1)<1$, the function
    \formula{
    \frac{1}{1-(ax+b+cx^{-1})G(x)}
    }
    is in $\mathscr{Q}$. We claim that $G_n^N\in\mathscr{Q}$ and $G_n^N(1)<1$ for every $y\leq n\leq N$.

    Our argument is by backward induction. We have $G_N^N(1) = g_N < 1$ and $G_N^N \in \mathscr{Q}$, since positive constant functions belong to $\mathscr{Q}$. For $\beta \in [0,1]$, we observe that
    \formula{
    T_\beta(x)= (1-\beta)^2 x + 2\beta(1-\beta) + \beta^2 x^{-1}.
    }
    Choosing
    \formula{
    & a = (1-g_n)(1-\beta_{n+1})^2, & b = (1-g_n)2\beta_{n+1}(1-\beta_{n+1}), && c = (1-g_n) \beta_{n+1}^2 &
    }
    yields
    \formula{
    & a+b+c = 1-g_n >0, & (a+b+c)G_{n+1}^N(1) = (1-g_n) G_{n+1}^N(1) < 1. &
    }
    It follows that $\smash{G^N_n/g_n}\in \mathscr{Q}$. The class $\mathscr{Q}$ is closed under multiplication by positive constants, so $\smash{G_n^N} \in \mathscr{Q}$. Moreover, by \eqref{eq:GN:recurrence}, $T_{\beta_{n+1}}(1)=1$ and $G_{n+1}^N(1)<1$, we have
    \formula{
    G_n^N(1) = \frac{g_n}{1-(1-g_n) T_{\beta_{n+1}}(1) G_{n+1}^N(1)} <\frac{g_n}{1-(1-g_n)} =1.
    }
    This completes the proof by induction.
\end{proof}

\begin{proposition}
\label{prop:G}
    For the family of functions $G_n^N$ defined by \eqref{eq:GN:recurrence}, we have
    \formula{
    \lim_{N \to \infty} G_y^N(x) = G_y(x)
    }
    uniformly on every compact subset of $\C \setminus[0,\infty)$.
\end{proposition}

\begin{proof} 
    The proof is divided into three steps.

    \textit{Step 1.} For $t>0$ and $w \in (0,1)$, denote
    \formula{
    & \Phi_n(w) = \frac{g_n}{1+c_n(t) w}, & c_n(t) = -(1-g_n) T_{\beta_{n+1}}(-t).&
    }
    Since $(G_n)$ is a family of nonconstant functions in $\G$, we have $g_n \in (0,1)$. Furthermore, for every $\beta \in [0,1]$,
    \formula{
    -T_{\beta}(-t) = (\beta t^{-1} - (1-\beta))(\beta - (1-\beta)t) = \frac{(\beta-(1-\beta)t)^2}{t}.
    }
    Using the inequality $|\beta - (1-\beta)t| \leq \max\{1,t\}$ for $\beta \in [0,1]$ and $t>0$, we deduce that
    \formula[eq:cn]{
    0 \leq c_n(t) \leq \max\{t, t^{-1}\}.
    }
    For $0 < a \leq b < 1$, we estimate
    \formula{
    |\log \Phi_n(b) - \log \Phi_n(a) | &= \int_a^b \frac{c_n}{1+c_n s} \, ds = \int_a^b \frac{c_ns}{1+c_n s} \, \frac{ds}{s} \\
    &\leq \int_a^b \frac{\max\{t, t^{-1}\}}{1+\max\{t, t^{-1}\}} \, \frac{ds}{s} = \frac{\max\{t,1\}}{1+t} \, |\log b - \log a|.
    }
    Above we used the monotonicity of $x \mapsto x/(1+x)$ and \eqref{eq:cn}. Setting $\rho_t = \max\{t,1\}/(1+t) \in (0,1)$, we conclude that for $0 < a \leq b < 1$,
    \formula{
    |\log \Phi_n(b) - \log \Phi_n(a) | \leq \rho_t |\log b - \log a|.
    }

    \textit{Step 2.} By Lemmas \ref{lem:G:bounds} and \ref{lem:GN:G}, we deduce that $0<G_n(-t)<1$ for $n \geq y$ and $0<G_n^N(-t)<1$ for $y \leq n \leq N$. Furthermore, by Corollary~\ref{cor:recurrence} and the definition of $G_n^N$, we have
    \formula[eq:Gns]{
    & G_n(-t)=\Phi_n(G_{n+1}(-t)), &G_n^N(-t)=\Phi_n(G_{n+1}^N(-t))
    }
    for $y \leq n < N$. Using \eqref{eq:Gns} and \eqref{eq:cn}, we find that for $y \leq n < N$,
    \formula{
    |\log G_n^N(-t) - \log G_n(-t) | &\leq \rho_t |\log G_{n+1}^N(-t) - \log G_{n+1}(-t) |.
    }
    Moreover, since $G_N^N(-t) = g_N$, \eqref{eq:Gn:recurrence} yields
    \formula{
    |\log G_N^N(-t) - \log G_N(-t) | &\leq \log \bigl(1+c_N(t) G_{N+1}(-t)\bigr) \\
     &\leq \log\bigl(1+\max\{t,t^{-1}\}\bigr).
    }
    By iteration, we arrive at
    \formula{
    |\log G_y^N(-t) - \log G_y(-t) | \leq \rho_t^{N-y} \log\bigl(1+\max\{t,t^{-1}\}\bigr).
    }
    Since $\rho_t<1$, the right-hand side tends to zero. Consequently, $G_y^N(-t) \to G_y(-t)$ as $N \to \infty$, for every $t>0$.

    \textit{Step 3.} We have already proved that the desired convergence holds on $(-\infty,0)$. We aim to extend this region by using Vitali's theorem (see Theorem 7.1.4 in \cite{remmert}): for a domain $D \subset \C$ and a sequence $(f_N)$ of holomorphic functions that is locally uniformly bounded in $D$, if the set
    \formula{
    \{x \in D: \lim_{N \to \infty} f_N(x) \text{ exists in } \C\}
    }
    has at least one accumulation point in $D$, then the sequence $(f_N)$ converges compactly in $D$.
    
    We set $D = \C \setminus[0,\infty)$ and $f_N = G_y^N$. By Lemma \ref{lem:GN:G}, we have $\smash{G_y^N} \in \G$, hence $f_N$ are holomorphic in $D$. In the previous step, we proved that the limit of $f_N$ exists in $(-\infty,0)$, which has accumulation points in $D$. It suffices to show that $f_N$ are locally uniformly bounded in $D$.

    Let $K$ be a compact subset of $D$. By the integral representation \eqref{eq:G:int:1} of $G_y^N$, there are measures $\smash{\mu_+^N}$ and $\smash{\mu_-^N}$ on $[0,1)$ satisfying $\smash{\mu_+^N[0,1) = \mu_-^N[0,1) = a^N(0) \leq 1}$, such that for $x \in K$,
    \formula{
    |G_y^N(x)| &= \bigg\lv \int_{[0,1)} \frac{1}{1-sx} \, \mu_+^N(ds) + \int_{[0,1)} \frac{1}{1-sx^{-1}} \, \mu_-^N(ds) - a^N(0) \bigg\rv \\
    &\leq \sup_{s \in [0,1)} \frac{1}{|1-sx|} a^N(0) + \sup_{s \in [0,1)} \frac{1}{|1-sx^{-1}|} a^N(0) + a^N(0) \\
    &\leq \sup_{s \in [0,1)} \frac{1}{|1-sx|} + \sup_{s \in [0,1)} \frac{1}{|1-sx^{-1}|} + 1.
    }
    The functions $(x,s)\mapsto |1-sx|^{-1}$ and $(x,s)\mapsto |1-sx^{-1}|^{-1}$ are continuous on the compact set $K\times[0,1]$, since their denominators do not vanish there. Consequently, the above expression is bounded by
    \formula{
    C_K = \sup_{(x,s)\in K\times[0,1]} \frac{1}{|1-sx|} + \sup_{(x,s)\in K\times[0,1]} \frac{1}{|1-sx^{-1}|} +1<\infty.
    }
    Thus $|G_y^N(x)|\leq C_K$ for every $N\geq y$ and $x\in K$. In particular, $C_K$ does not depend on $N$. It follows that
    \formula{
    \sup_{N \geq y} \sup_{x \in K} |G_y^N(x)| < \infty
    }
    on every compact $K \subset D$.
    
    Vitali's theorem implies that there is a function $\tilde G_y$ holomorphic in $D$ such that $\smash{G_y^N \to \tilde G_y}$ uniformly on every compact subset. On the other hand, we have seen that $\smash{G_y^N \to G_y}$ on $(-\infty,0)$. Since $(-\infty,0)$ has accumulation points in $D$, we deduce that $G_y = \tilde G_y$ on $D$.
\end{proof}

\section{Proofs of the main results}
\label{sec:proofs}
\subsection{Proof of Theorem \ref{thm:main}}

\begin{proof}[Proof of Theorem \ref{thm:main} (sufficient condition)]
Define $y=2A+B+C+1$ and let $F$ be a function given by \eqref{eq:main:F}. The idea of the proof comes from \cite{bondesson}, where the proof of an analogous one-dimensional theorem is conducted in two steps (see Theorem 2 therein). We show that there is a random walk such that
\formula[]{
F_y^{y-1,\infty}(x) &= G(x), \label{eq:F1:proof} \\
F_{y-1}^{0,y}(x) &= \prod_{k=1}^A \frac{1-q_k}{1-q_k x} \prod_{k=1}^B \frac{1-r_k}{1-r_k/x} \prod_{k=1}^{C} (\beta_k x^{-1}+1-\beta_k).\label{eq:F2:proof}
}
If we put $p_{y-1,1}=p_{y-1,2}=0$, then the result will follow from \eqref{eq:F:factors}.


\textit{Step 1.} We start from the difficult step, which is proving \eqref{eq:F1:proof}. Applying Corollary \ref{cor:recurrence} to $G_y=G$, we find that there are functions $G_n \in \G \cup \{0\}$ with $g_n = \sup_{x<0}G_n(x)$, and constants $\beta_{n+1} \in [0,1]$, such that for every $n \geq y$,
\formula[eq:Gn:proof]{
    G_n(x) = \frac{g_n}{1-(1-g_n)T_{\beta_{n+1}}(x) G_{n+1}(x)},
}
where
\formula{
    T_\beta(x) = (\beta x^{-1} + 1-\beta)((1-\beta)x+\beta).
}
If we denote $H_n(x) = (\beta_n x^{-1} + 1-\beta_n) G_n(x)$, then \eqref{eq:Gn:proof} can be rewritten as
\formula[eq:Hn:recursion]{
H_n(x) = \frac{g_n(\beta_n x^{-1} + 1-\beta_n)}{1-(1-g_n)((1-\beta_{n+1})x + \beta_{n+1}) H_{n+1}(x)},
}
where $\beta_y=0$. Comparing the above relation with \eqref{eq:F:b}, we deduce that for $n \geq y$ the functions $H_n$ and $F_n^{n-1,\infty}$ satisfy the same recurrence when the transition probabilities are chosen as follows:
\formula*[eq:probs]{
    & p_{n,1} = (1-g_n)(1-\beta_{n+1}), \quad &&p_{n,2} = (1-g_n)\beta_{n+1}, \\ 
    & p_{n,3} = g_n \beta_n, &&p_{n,4} = g_n(1-\beta_n).
}

Suppose that $G_N \equiv g_N \in(0,1]$ for some $N\geq y$. By our convention, $G_n\equiv0$ and $g_n=0$ for every $n>N$. Hence $p_{n,3}=p_{n,4}=0$ for all $n>N$, so the walk cannot return to level $N$ after reaching level $N+1$. Consequently, a walk starting at level $N$ reaches level $N-1$ if and only if its first step is downwards. Thus,
\formula{
F_N^{N-1,\infty}(x) = p_{N,3}x^{-1}+p_{N,4} = g_N(\beta_Nx^{-1}+1-\beta_N) = H_N(x).
}
Using the common recurrence, backward induction yields $F_n^{n-1,\infty}=H_n$ for every $y\leq n\leq N$. In particular, \eqref{eq:F1:proof} holds, and both sides are rational functions.

Assume now that $G_n$ are nonconstant for every $n\geq y$. In this case, the recurrence does not terminate. Without any terminal condition the solutions need not be unique. To resolve this issue, we consider the modified random walk with transition probabilities given by
\formula{
p_{k,j}^{N}=p_{k,j},
}
for $y\leq k\leq N$, $j\in\{1,2,3,4\}$, where $p_{k,j}$ are chosen as in \eqref{eq:probs}, and
\formula{
& p_{k,1}^{N}=1, & p_{k,2}^{N}=p_{k,3}^{N}=p_{k,4}^{N}=0 &
}
for $k \geq N+1$. Thus, once the walk reaches level $N+1$, it moves upwards at every subsequent step and never returns to level $N$. For first passage events below level $N+1$, this construction is equivalent to killing the original walk upon reaching that level.

For $y\leq n\leq N$, let $\smash{F_n^{N}}$ denote the generating function of the first passage location from level $n$ to level $n-1$ for this walk. Starting from level $N$, the walk reaches level $N-1$ if and only if its first step is downwards. Hence,
\formula{
F_N^{N}(x) = p_{N,3}x^{-1}+p_{N,4} = g_N (\beta_N x^{-1}+1-\beta_N).
}
Let $G_n^N$ be the functions defined in Proposition~\ref{prop:G}. The terminal identity above and the recurrence \eqref{eq:F:b} imply, by backward induction, that
\formula{
F_n^N(x)=(\beta_nx^{-1}+1-\beta_n)G_n^N(x)
}
for $y\leq n\leq N$.

On the one hand, by Proposition~\ref{prop:G}, we have $G_y^N(x) \to G_y(x)=G(x)$, as $N \to \infty$, uniformly on every compact subset of $\C \setminus[0,\infty)$. In particular, for every $|x|=1$ and $x \neq 1$,
\formula{
\lim_{N \to \infty} F_y^N(x) = G(x).
}
On the other hand, as we have already seen in \eqref{eq:F:limit}, for $|x|=1$,
\formula{
\lim_{N \to \infty} F_y^N(x) = F_y^{y-1,\infty}(x).
}
Therefore,
\formula[eq:equal]{
F_y^{y-1,\infty}(x) = G(x)
}
for $|x|=1$, $x\neq1$. Both sides are continuous on the unit circle, since the corresponding sequences are summable. Hence, \eqref{eq:equal} extends to the whole unit circle $|x|=1$, which completes the proof of \eqref{eq:F1:proof}.

\textit{Step 2.} We now prove \eqref{eq:F2:proof}. Recall that in the first step of the proof, we chose transition probabilities $p_{n,k}$ for $n \geq y$. Furthermore, we set $p_{y-1,1}=p_{y-1,2}=0$, which means that once the random walk reaches level $y-1$, it never goes back above. Our construction involves three types of blocks.

The first type of block realises geometric distributions. For $n \in \N$ and $q \in (0,1)$, we define
\formula{
& p_{n,4}=1, & p_{n-1,1}=q, && p_{n-1,4}=1-q &
}
and the remaining probabilities on levels $n$, $n-1$ are equal to zero. Then the corresponding generating function is given by
\formula{
F_{n}^{n-2,\infty}(x) = \frac{1-q}{1-qx}.
}
We build $A$ blocks of this type. Each block consists of two levels and in each block we use one of the parameters $q_1, q_2, \ldots, q_A$. Since $y-1=2A+B+C$, this gives
\formula[eq:F1]{
F_{2A+B+C}^{B+C,\infty}(x) = \prod_{k=1}^A \frac{1-q_k}{1-q_kx}.
}

The second type of block realises reflected geometric distributions. However, a reflected geometric distribution alone cannot be realised with our two-level construction with a barrier. To solve this problem, we pair each reflected geometric factor with a reflected two-point factor. For $n \in \N$, $r \in (0,1)$ and $\beta \in (0,1]$, we define
\formula{
& p_{n,3}=\beta, && p_{n,4}=1-\beta, && p_{n-1,2} = \frac{r}{\beta+r(1-\beta)}, & p_{n-1,4} = \frac{\beta(1-r)}{\beta+r(1-\beta)} &
}
and again the remaining probabilities on these levels are set to zero. The corresponding generating function reads
\formula{
F_n^{n-2,\infty}(x) = \frac{1-r}{1-rx^{-1}} (\beta x^{-1}+1-\beta).
}
We build $B$ blocks of this type. Each block uses two levels and one parameter from each of the lists $r_1, r_2, \ldots, r_B$ and $\beta_1, \beta_2, \ldots, \beta_B$. In particular, here we use the assumption that $B \leq C$, which provides enough parameters $\beta$ to do this. Consequently,
\formula[eq:F2]{
F_{B+C}^{C-B,\infty}(x) = \prod_{k=1}^B \frac{1-r_k}{1-r_kx^{-1}} \prod_{k=1}^B (\beta_k x^{-1}+1-\beta_k).
}

If $C=B$, our construction is finished. Otherwise, we construct $C-B$ one-level blocks. For $n \in \N$, we define
\formula{
& p_{n,3}=\beta, & p_{n,4} =1-\beta &
}
with the rest of the probabilities vanishing. On each level, we use one of the remaining parameters $\beta_{B+1}, \beta_{B+2}, \ldots, \beta_C$. Then the corresponding generating function is
\formula[eq:F3]{
F_{C-B}^{0,\infty}(x) = \prod_{k=B+1}^{C} (\beta_k x^{-1}+1-\beta_k).
}

By the strong Markov property and the barriers separating the blocks,
\formula{
F_{y-1}^{0,y}(x) &= F_{2A+B+C}^{B+C,\infty}(x) \, F_{B+C}^{C-B,\infty}(x) \, F_{C-B}^{0,\infty}(x).
}
Our assertion now follows from \eqref{eq:F1}, \eqref{eq:F2} and \eqref{eq:F3}.
\end{proof}

It remains to prove the converse.

\begin{proof}[Proof of Theorem \ref{thm:main} (necessary condition)]
Consider a diagonal random walk starting from $(0,y)$ for some $y \in \N$ with a nonzero generating function $F$ of the first passage location. By Theorem 1.1 in \cite{fpl}, $F$ can be represented as
\formula[eq:F:proof2]{
F(x)=G(x)\prod_{k=1}^{y}\frac{p_{k,3}x^{-1}+p_{k,4}}{p_{k,3}+p_{k,4}}
\prod_{k=1}^{A_y}\frac{1-q_k}{1-q_kx}
\prod_{k=1}^{B_y}\frac{1-r_k}{1-r_kx^{-1}},
}
where $G\in\mathscr{G}$, $q_k,r_k\in(0,1)$, and $A_y$, $B_y$ are nonnegative integers. Furthermore, $p_{k,3}+p_{k,4}>0$ for every $1 \leq k \leq y$; otherwise $F=0$. We set
\formula{
\beta_k=\frac{p_{k,3}}{p_{k,3}+p_{k,4}}
}
and observe that $\beta_k=0$ if and only if $p_{k,3}=0$, and $\beta_k = 1$ if and only if $p_{k,4}=0$. Consequently, $F$ can be represented as \eqref{eq:main:F} with $A=A_y$, $B=B_y$ and $C=\sum_{k=1}^y \ind_{p_{k,3}>0}$. It remains to prove that
\formula[eq:B:bound]{
B_y \leq \sum_{k=1}^y \ind_{p_{k,3}>0}.
}

For a nonzero Laurent polynomial $H$, let $d_+(H)$ denote its highest exponent of $x$, and let $d_-(H)$ denote the negative of its lowest exponent. For $0\leq n\leq y$, define
\formula{
H_n(x)=\frac{1}{F_n^{0,n+1}(x)}\prod_{k=1}^{n}(p_{k,3}x^{-1}+p_{k,4}).
}
By Lemma~5.1 in \cite{fpl}, with the barrier cases treated in Section~5.2 therein, these functions are Laurent polynomials admitting the factorisation
\formula{
H_n(x)=H_n(1)
\prod_{k=1}^{d_+(H_n)}\frac{1-q_{n,k}x}{1-q_{n,k}}
\prod_{k=1}^{d_-(H_n)}\frac{1-r_{n,k}x^{-1}}{1-r_{n,k}},
}
for appropriate sequences $q_{n,k},r_{n,k}\in(0,1)$ and $H_n(1)>0$. The geometric factors in \eqref{eq:F:proof2} arise from $H_{y-1}(1)/H_{y-1}(x)$, so $A_y=d_+(H_{y-1})$ and $B_y=d_-(H_{y-1})$.

Furthermore, $H_0=H_1=1$, and Eq.~(5.3) in \cite{fpl} gives for $2 \leq n \leq y$,
\formula[eq:Hs]{
H_n(x)=H_{n-1}(x)-(p_{n-1,1}x+p_{n-1,2})(p_{n,3}x^{-1}+p_{n,4})H_{n-2}(x).
}
The coefficient multiplying $H_{n-2}$ contains a term in $x$ only if $p_{n,4}>0$, and a term in $x^{-1}$ only if $p_{n,3}>0$. Consequently,
\formula{
d_+(H_n)&\leq\max\big\{d_+(H_{n-1}),d_+(H_{n-2})+\ind_{p_{n,4}>0}\big\},\\
d_-(H_n)&\leq\max\big\{d_-(H_{n-1}),d_-(H_{n-2})+\ind_{p_{n,3}>0}\big\}.
}
Since $d_\pm(H_0)=d_\pm(H_1)=0$, induction yields
\formula[deg:bounds]{
& d_+(H_n)\leq\sum_{k=2}^{n}\ind_{p_{k,4}>0}, & d_-(H_n)\leq\sum_{k=2}^{n}\ind_{p_{k,3}>0}. &
}
Now \eqref{eq:B:bound} follows from taking $n=y-1$.
\end{proof}

\subsection{Proof of Theorem~\ref{thm:main:2}}

\begin{proof}[Proof of Theorem~\ref{thm:main:2}]

Let $F$ be a nonzero generating function of first passage location. Then, by Theorem \ref{thm:main}, $F$ can be represented as \eqref{eq:main:F} with $B \leq C$. Assume that a maximal representation of $F$ is
\formula{
F(x)=\tilde G(x)
\prod_{k=1}^{\tilde A}\frac{1-\tilde q_k}{1-\tilde q_kx}
\prod_{k=1}^{\tilde B}\frac{1-\tilde r_k}{1-\tilde r_k x^{-1}} \prod_{k=1}^{\tilde{C}} (\tilde\beta_kx^{-1} + 1-\tilde\beta_k ),
}
where $\tilde G\in\mathscr{G}$, $\tilde q_k, \tilde r_k\in(0,1)$, $\tilde \beta_k \in (0,1]$ and $\tilde A, \tilde B, \tilde C$ are nonnegative integers. Moreover, the functions
\formula{
& \tilde G(x)\frac{1-\tilde q_k}{1-\tilde q_kx}, & \tilde G(x) \frac{1-\tilde r_k}{1-\tilde r_kx^{-1}} &
}
are not in $\G$ for any $1 \leq k \leq \tilde A$ and $1 \leq k \leq \tilde B$, respectively.

\textit{Step 1.} 
We show that $C=\tilde C$.

Since $G$, $\tilde G$ and all geometric factors are strictly positive on $(-\infty,0)$, the negative zeros of $F$, counted with multiplicity, determine the parameters $\beta_k\in(0,1)$ through the relation $x=-\beta_k/(1-\beta_k)$. Hence these parameters coincide in both representations. Cancel the corresponding factors.

By Lemma \ref{lem:reG}, $\re G(x)>0$ and $\re \tilde G(x)>0$ for $|x|=1$. All geometric factors also have positive real parts on the unit circle. Consequently, each of these factors has zero change of argument along the unit circle, whereas each remaining two-point factor, with parameter $\beta=1$, equals $x^{-1}$ and contributes $-2\pi$. Equality of the remaining products implies that the numbers of parameters $\beta$ equal to $1$ coincide. Thus $C=\tilde C$.

\textit{Step 2.} We prove that $\tilde B \leq B$.

Assume that $\tilde B > B$ and observe that
\formula[eq:Gs]{
\tilde G(x) \prod_{k=1}^{\tilde A} \frac{1-\tilde q_k}{1-\tilde q_k x} \prod_{k=1}^{\tilde B} \frac{1-\tilde r_k}{1-\tilde r_k x^{-1}} = G(x) \prod_{k=1}^{A} \frac{1-q_k}{1-q_k x} \prod_{k=1}^{B} \frac{1-r_k}{1-r_k x^{-1}}.
}
Using the formulas
\formula{
\log(1-qx) = -\int_{1/q}^\infty \biggl( \frac{1}{s-x} - \frac{s}{s^2+1}\biggr) \, ds + \log \sqrt{q^2+1}, \\
\log(1-rx^{-1}) = \int_{0}^r \biggl( \frac{1}{s-x} - \frac{s}{s^2+1}\biggr) \, ds + \log \sqrt{r^2+1}
}
and exponential representations \eqref{eq:G:exp} of $G$ and $\tilde G$ with corresponding functions $\ph$ and $\tilde \ph$, respectively, we express both sides of \eqref{eq:Gs} in exponential form. Since the corresponding exponential densities are equal almost everywhere, we have
\formula{
-\tilde\ph(s) - \sum_{k=1}^{\tilde A} \ind_{(1/\tilde q_k, \infty)}(s) + \sum_{k=1}^{\tilde B} \ind_{(0, \tilde r_k)}(s) = -\ph(s)- \sum_{k=1}^{A} \ind_{(1/q_k,\infty)}(s) + \sum_{k=1}^{B} \ind_{(0,r_k)}(s).
}
    Let $\tilde r = \min\{\tilde r_k: 1 \leq k \leq \tilde B \}$. Recall that $-1 \leq \ph, \tilde\ph \leq 0$ on $(0,1)$. Moreover, since $B, \tilde B$ are integers, our assumption that $\tilde B > B$ implies $\tilde B \geq B+1$. Taking $s \in (0,\tilde r)$ yields, for almost every $s$,
\formula{
\tilde B \leq - \tilde\ph(s) + \tilde B \leq -\ph(s) + B \leq 1+ B \leq \tilde B.
}
It follows that $\tilde \ph(s)=0$ almost everywhere for $s \in (0,\tilde r)$. Hence, the function
\formula{
H(x)=\tilde G(x) \, \frac{1-\tilde r}{1-\tilde rx^{-1}}
}
admits exponential representation \eqref{eq:G:exp}. Its exponential density is $\tilde\ph(s)-\ind_{(0,\tilde r)}(s)$. Since $\tilde\ph=0$ almost everywhere on $(0,\tilde r)$, this takes values in $[-1,0]$ on $(0,1)$ and in $[0,1]$ on $(1,\infty)$. Furthermore, $H$ is continuous on the unit circle. Hence,
\formula{
\lim_{t\to0}H(e^{it})=H(1)=\tilde G(1)\in(0,1].
}
Together with the exponential representation obtained above, Theorem~\ref{thm:G} implies that $H\in\G$, contradicting maximality of the representation..

\textit{Step 3.}
The previous steps and the assumption that $B \leq C$ imply that
\formula{
\tilde B \leq B \leq C = \tilde C,
}
as desired.
\end{proof}

\section{Discussion and examples}
\label{sec:discussion}

In \cite{fpl}, other types of random walks were considered. The distribution of the first passage location was found for standard random walk (moving up, down, right or left) and random walk on the honeycomb lattice, see Theorem 1.2 and Corollary 1.3 therein, respectively.

The converse of Theorem~1.2 in \cite{fpl} for standard random walks requires additional restrictions on the $\amcm$ factor; see Example \ref{ex:news:not}. The geometric part alone, however, can be realised by a simple choice of transition probabilities as in Example \ref{ex:news:geo}. No two-point factors occur in the corresponding representation.

\begin{example}
\label{ex:news:not}
Consider
\formula{
F(x)=\frac{1}{\sqrt{2-x}}.
}
Note that $F$ admits the representation \eqref{eq:main:F} with $A=B=C=0$ and $F=G$. Indeed,
\formula{
G(x) = \frac{1}{\pi\sqrt{2}}\int_0^1\frac{2}{2-sx}\,\frac{ds}{\sqrt{s(1-s)}}
}
is in $\G$ by \eqref{eq:G:int:1}. This representation is maximal. However, $F$ cannot be realised as a generating function of first passage location for a standard random walk from any starting level $y$.

Suppose otherwise and assume that there is a standard random walk $(X_n, Y_n)$ starting from $(0,y)$ such that $F_y^{0,\infty} = F$, where for $b \leq a$ we define
\formula{
F_a^{b,\infty}(x) = \E^{(0,a)} \bigl[ x^{X(\tau_b)} \ind_{\tau_b < \infty}\bigr].
}
If $a=n$ and $b=n-1$, the corresponding functions are in $\G$ for every $n \leq y$ (see Section 6.1 in \cite{fpl}). Denote $\smash{H_n = F_n^{n-1,\infty}}$. In particular, $\smash{0 < H_n(-t)\leq 1}$ for $t>0$. By the strong Markov property, one can also deduce that
\formula[eq:F:prod]{
F_y^{0,\infty}(x) = \prod_{n=1}^y H_n(x).
}
Moreover, by Eq. (6.1) in Section 6.1 in \cite{fpl}, we have
\formula[eq:F:rec]{
H_n(x)
=\frac{p_{n,4}}{1-p_{n,1}x-p_{n,3}x^{-1}-p_{n,2}H_{n+1}(x)},
}
where the transition probabilities are numbered E, N, W, S. Necessarily, $p_{n,4}>0$ for $1\leq n\leq y$.

If $p_{n,1}=0$ for every $1\leq n\leq y$, then using  \eqref{eq:F:prod} and \eqref{eq:F:rec}, we get for $t>0$,
\formula{
F(-t) = \prod_{n=1}^y H_n(-t) = \prod_{n=1}^y \frac{p_{n,4}}{1+p_{n,1}t+p_{n,3}t^{-1}-p_{n,2}H_{n+1}(-t)} \geq \prod_{n=1}^y \frac{p_{n,4}}{1+p_{n,3}t^{-1}}.
}
Observe that
\formula{
& \lim_{t\to\infty} F(-t) = \lim_{t\to\infty} \frac{1}{\sqrt{2+t}}=0, & \lim_{t\to\infty} \prod_{n=1}^y \frac{p_{n,4}}{1+p_{n,3}t^{-1}} = \prod_{n=1}^y p_{n,4}>0.
}
This is a contradiction. Thus, $p_{n,1}>0$ for some $1\leq n\leq y$. For this index, \eqref{eq:F:prod} and \eqref{eq:F:rec} imply that
\formula{
F(-t)\leq H_n(-t)= \frac{p_{n,4}}{1+p_{n,1}t+p_{n,3}t^{-1} - p_{n,2}H_{n+1}(-t)} \leq \frac{p_{n,4}}{p_{n,1}t}.
}
Since
\formula{
\lim_{t \to \infty} tF(-t) = \lim_{t \to\infty} \frac{t}{\sqrt{2+t}} = \infty,
}
we again have a contradiction. No standard random walk from any starting level $y$ satisfies $F_y^{0,\infty} = F$.
\end{example}

\begin{example}
\label{ex:news:geo}
Consider
\formula{
F(x) = \prod_{k=1}^A \frac{1-q_k}{1-q_kx} \prod_{k=1}^B \frac{1-r_k}{1-r_kx^{-1}}
}
for nonnegative integers $A,B$ and $q_k, r_k \in (0,1)$. We will construct a standard random walk such that $F_y^{0,\infty}=F$ for a suitable $y$.

Let $y=A+B$. As in the diagonal case, we construct two types of blocks. Using the notation from Example \ref{ex:news:not}, we define for $1 \leq k \leq A$,
\formula{
& p_{y-k+1,1} = q_k, & p_{y-k+1,4} = 1-q_k, &
}
and for $1 \leq k \leq B$,
\formula{
& p_{k,3} = r_k, & p_{k,4} = 1-r_k. &
}
All remaining transition probabilities are set to zero. In other words, in the first $A$ blocks, the random walk can jump either right or down. In the remaining $B$ blocks, we allow moving either left or down. The corresponding generating functions are
\formula{
& F_{A+B}^{B,\infty}(x) = \prod_{k=1}^A \frac{1-q_k}{1-q_kx}, &  F_{B}^{0,\infty}(x) = \prod_{k=1}^B \frac{1-r_k}{1-r_kx^{-1}}. &
}
Therefore,
\formula{
F_{A+B}^{0, \infty}(x) = F_{A+B}^{B,\infty}(x) F_{B}^{0,\infty}(x)=F(x).
}
\end{example}

Under the embedding described in Section~6.2 of \cite{fpl}, random walks on the honeycomb lattice can be represented as diagonal random walks with certain transition probabilities fixed at zero. A realisation on the honeycomb lattice must also respect these constraints. The following example shows that they can exclude a generating function covered by Theorem~\ref{thm:main}.

\begin{example}
\label{ex:honeycomb:not}
Consider the function
\formula{
F(x)=\frac{1}{\sqrt{2-x}}
}
from Example~\ref{ex:news:not}. It is the generating function of an $\amcm$ sequence and satisfies $F(1)=1$. Nevertheless, it cannot be realised as the generating function of the rescaled first passage location of a honeycomb random walk, regardless of the starting level.

We use the diagonal embedding from Section 6.2 of \cite{fpl}, in which
\formula{
p_{n,1}=0 \quad \text{for odd }n,
\qquad
p_{n,3}=0 \quad \text{for even }n.
}
Suppose that $F_y^{0,\infty}=F$ for some $y\geq1$, and write $H_n=F_n^{n-1,\infty}$. We have $p_{n,3}+p_{n,4}>0$ for all $1\leq n\leq y$.

For even $n$, the downward two-point factor is constant, so $H_n\in\G \cup\{0\}$. In particular, $0\leq H_n(-t)\leq1$ for $t>0$, and $H_n(-t)$ has a finite limit as $t\to\infty$. The recurrence \eqref{eq:F:b} gives
\formula{
H_n(-t)=
\begin{cases}
\displaystyle \frac{p_{n,4}-p_{n,3}t^{-1}}{1-p_{n,2}H_{n+1}(-t)}, & \text{if } n \text{ is odd},\\[8pt]
\displaystyle \frac{p_{n,4}}{1+(p_{n,1}t-p_{n,2})H_{n+1}(-t)}, & \text{if }n \text{ is even}.
\end{cases}
}

For odd $n$ with $H_n\not\equiv0$, the denominator has a finite limit and is bounded below by
\formula{
1-p_{n,2}=p_{n,3}+p_{n,4}>0.
}
Consequently, if $p_{n,4}>0$, then $H_n(-t)$ converges to a number in $(0,1]$. If $p_{n,4}=0$, then $t H_n(-t)$ converges to a number in $[-1,0)$. The remaining possibility is $H_n\equiv0$.

Now let $n\leq y$ be even, so that $p_{n,4}>0$. If $p_{n,1}>0$ and $H_{n+1}(-t)$ has a positive limit, the denominator in the second recurrence is asymptotic to a positive constant times $t$. In every other case, it has a finite positive limit. Indeed, when $p_{n,1}=0$, this limit is at least $1-p_{n,2}\geq p_{n,4}$; when $H_{n+1}(-t)\to0$, it is at least $1-p_{n,1}\geq p_{n,4}$.

Thus, for every $1\leq n\leq y$, the function $H_n(-t)$ is asymptotic either to a nonzero constant or to a nonzero constant times $t^{-1}$. By the strong Markov property, it follows that
\formula{
F(-t)=F_y^{0,\infty}(-t)=\prod_{n=1}^y H_n(-t)
\approx t^{-m}
}
for an integer $m\in\{0,\ldots,y\}$. This contradicts $F(-t)\approx t^{-1/2}$. Hence no such honeycomb random walk exists.
\end{example}

The examples above show that the converse problem for standard and honeycomb random walks requires additional restrictions on the $\amcm$ component. A natural direction for further research is to identify necessary and sufficient conditions on this component. The necessary asymptotic restrictions derived above provide a starting point for such a characterisation.

\section*{Acknowledgements}
OpenAI's Codex was used to assist with reviewing the mathematical arguments, checking references, revising parts of the manuscript, developing counterexamples and exploring proof strategies for Proposition~\ref{prop:G}.

\bibliographystyle{abbrv}
\bibliography{bibi}

\end{document}